\documentclass[]{article}    % 939-46-431 = Abstract Number (Sep 24, 1999)
 \usepackage{spie}
 \usepackage{amsfonts}
 \usepackage[mathscr]{eucal}

\newtheorem{corollary}[theorem]{Corollary}
\newtheorem{example}[theorem]{Example}
\newtheorem{proposition}[theorem]{Proposition}
\newtheorem{problem}[theorem]{Problem}

\begin{document}
%%%%%%%%%%%%%%%%%%%%%%%%%%%%%%%%%%%%%%%%%%%%%%%%%%%%%%%%%%%%%%%%%%%%%%%%%%%%%
\title{Modular frames for Hilbert C*-modules and symmetric approximation
  of frames}

\author{Michael Frank\supit{a} and David R.~Larson\supit{b} 
\skiplinehalf 
\supit{a}Universit\"at Leipzig, Mathematical Institut, D-04109 Leipzig, F.R.~Germany
\\
\supit{b}Dept.~Mathematics, Texas A{\&}M University, College Station, TX 77843, U.S.A.
}

%\author{Michael Frank \\ [12pt]
%  Universit\"at Leipzig, Math.~Institut \\ D-04109 Leipzig, F.R.~Germany \\
%  [12pt]
%  David R.~Larson \\ [12pt]
%  Dept.~Mathematics, Texas A{\&}M Univ. \\ College Station, TX 77843, U.S.A.}

\authorinfo{E-mails: frank@mathematik.uni-leipzig.de, larson@math.tamu.edu.\\
   \hspace*{0.5cm}
   Supported in part by grants from the Deutsche Forschungsgemeinschaft and
   the National Science Foundation.}
\pagestyle{plain}
%%%%%%%%%%%%%%%%%%%%%%%%%%%%%%%%%%%%%%%%%%%%%%%%%%%%%%%%%%%%%%%%%%%%%%%%%%%%%
\maketitle
%%%%%%%%%%%%%%%%%%%%%%%%%%%%%%%%%%%%%%%%%%%%%%%%%%%%%%%%%%%%%%%%%%%%%%%%%%%%%
\begin{abstract}
  We give a comprehensive introduction to a general modular frame
  construction in Hilbert C*-modules and to related linear operators on them.
  The Hilbert space situation appears as a special case. The reported
  investigations rely on the idea of geometric dilation to standard Hilbert
  C*-modules over unital C*-algebras that admit an orthonormal modular Riesz
  basis. Interrelations and applications to classical frame theory are
  indicated. Resorting to frames in Hilbert spaces we discuss some measures for
  pairs of frames to be close to one another. In particular, the existence and
  uniqueness of the closest (normalized) tight frame to a given frame is
  investigated. For Riesz bases with certain restrictions the set of closest
  tight frames often contains a multiple of its symmetric orthogonalization.
\end{abstract}

%%%%%%%%%%%%%%%%%%%%%%%%%%%%%%%%%%%%%%%%%%%%%%%%%%%%%%%%%%%%%%%%%%%%%%%%%%%%%

\section{INTRODUCTION}

\noindent
The inner structure of Hilbert spaces is easily described by fixing a basis,
orthonormalizing it and working with the coordinates of every element with
respect to the latter. Considering finite-dimensional Hilbert spaces as free
$\mathbb C$-modules one could ask whether similar generating sets of finitely
generated projective C*-modules can be indicated which characterize the
modules up to isomorphism. Unfortunately, there are two obstacles: projective
C*-modules need not to be free in general, and there does not exist any
general notion of `C*-linear independence' of sets of generators because
of the existence of zero-divisors in any non-trivial C*-algebra. Beside
these circumstances one does not know any canonical method to replace the
process of Gram-Schmidt or symmetric orthogonalization of bases in the
situation of sets of modular generators.

In 1997, Working with frames for Hilbert spaces that arise canonically in wavelet
and Weyl-Heisenberg / Gabor frame theory, we got the idea to investigate
modular frames of Hilbert C*-modules over unital C*-algebras as a possible
replacement for the questionable analogs of bases. This class of
C*-modules includes finitely generated projective ones.
The resulting theory\cite{FL:98,FL:00} has been encouraging because of its
consistency and strength, and also because of the number of mathematical
problems which can make use of it. In the field of wavelet and frame theory
and its applications to signal and image processing many frames arise as
the result of group actions on single functions. Extending the group to
its (reduced) group C*-algebra or to its group von Neumann algebra and
taking the generated frame as a generating set of a Hilbert C*-module over
one of these C*-algebras we are in the context in which our concept can be
applied. This point of view has been of interest e.g.~to M.~A.~Rieffel, to
O.~Bratteli and P.~E.~T.~J{\o}rgensen, and to P.~G.~Casazza and M.~Lammers
as we know from ongoing discussions. In the literature we found other fields
of applications like the description of conditional expectations of finite
(Jones) index on C*-algebras (\cite{Wata,PP,BDH,FK}), the analysis of
Cuntz-Krieger-Pimsner algebras (\cite{DPZ,KPW}), the investigation of the
stable rank of C*-algebras (\cite{Badea,Vill}) and the search for $L^2$-invariants
in global analysis (\cite{CM,Lue,Fr}).

The purpose of the present paper is to give a survey on our results on
modular frames for Hilbert C*-modules indicating their generality and
strength as well as pointing out differences to the Hilbert space frame
theory and open problems. For full proofs and more details we refer to
our basic publications (\cite{FL:98,FL:00}). The next section we explain
decomposition and reconstruction results. In the third section we deal with
frame-related invariants of finitely generated projective C*-modules that
characterize them up to isomorphism. the fourth section is devoted to a
structure theorem on the nature of operators $\{ b_i \}_i$ on a certain
Hilbert space such that $\sum_i b_i^*b_i = id$. The last section is concerned
with a discussion on various problems of frame approximation by (normalized)
tight ones.
%%%%%%%%%%%%%%%%%%%%%%%%%%%%%%%%%%%%%%%%%%%%%%%%%%%%%%%%%%%%%%%%%%%%%%%%%%%%%

\section{MODULAR FRAMES FOR HILBERT C*-MODULES}

\noindent
The concept of Hilbert C*-modules arose as a generalization of the notions
`Hilbert space', `fibre bundle' and `ideal'. The basic idea has been to
consider modules over arbitrary C*-algebras instead of linear spaces and to
allow the inner products to take values in those C*-algebras of coefficients
being C*-(anti-)linear in their arguments.
For the history and for comprehensive accounts we refer to the publications
by E.~C.~Lance\cite{Lance}, by N.-E.~Wegge-Olsen\cite{NEWO} and by I.~Raeburn,
D.~P.~Williams\cite{RW}.

\begin{definition}
 Let $A$ be a (unital) C*-algebra and $\mathcal M$ be a (left) $A$-module.
 Suppose that the linear structures given on $A$ and $\mathcal M$ are
 compatible, i.e. $\lambda (a x) =(\lambda a)x = a(\lambda x)$ for every
 $\lambda \in {\mathbb C}$, $a \in A$ and $x \in \mathcal M$. If there
 exists a mapping $\langle .,. \rangle: {\mathcal M} \times {\mathcal M}
 \rightarrow A$ with the properties
   \newcounter{marke}
   \begin{list}{(\roman{marke})}{\usecounter{marke}}
      \item $\langle x,x \rangle \geq 0$ for every $x \in \mathcal M$,
      \item $\langle x,x \rangle =0$ if and only if $x=0$,
      \item $\langle x,y \rangle = \langle y,x \rangle^*$ for every $x,y \in
            \mathcal M$,
      \item $\langle ax,y \rangle = a \langle x,y \rangle$ for every $a \in A$,
            every $ x,y \in \mathcal M$,
      \item $\langle x+y,z \rangle = \langle x,z \rangle + \langle y,z \rangle$
            for every $x,y,z \in \mathcal M$,
   \end{list}
 then the pair $\{ \mathcal M, \langle .,. \rangle \}$ is called a
 {\it (left) pre-Hilbert $A$-module}. The map $\langle .,. \rangle$ is said to
 be an {\it $A$-valued inner product}.    
 If the pre-Hilbert $A$-module $\{ \mathcal M, \langle .,. \rangle \}$ is
 complete with respect to the norm $\|x\| = \| \langle x,x \rangle \|^{1/2}$
 then it is called a {\it Hilbert $A$-module}.

 Two Hilbert $A$-modules are {\it unitarily isomorphic} if there exists a
 bounded $A$-linear isomorphism of them which preserves the inner product
 values.

 In case $A$ is unital the Hilbert $A$-module $\mathcal M$ is {\it
 (algebraically) finitely generated} if there exists a finite set $\{ x_i
 \}_{i \in \mathbb N} \subset \mathcal M$ such that $x = \sum_i a_i x_i$ for
 every $x \in \mathcal M$ and some coefficients $\{ a_i \} \subset A$.
 If $A$ is unital the Hilbert $A$-module $\mathcal M$ is {\it countably
 generated} if there exists a countable set $\{ x_i \}_{i \in \mathbb N}
 \subset \mathcal M$ such that the set of all finite $A$-linear combinations
 $\{ \sum _j a_jx_j \}$, $\{ a_i \} \subset A$, is norm-dense in $\mathcal M$.
\end{definition}

Hilbert C*-modules appear naturally in a number of situations. For example,
the set of all essentially bounded measurable maps of a measure space $X$ into
a measurable field of Hilbert spaces on $X$ becomes a Hilbert
$L^\infty(X)$-module after factorization by the set of maps that are non-zero
only on sets of measure zero, (\cite{Fr92}). A simple algebraic construction
is the algebraic tensor product $A \odot H$ of a C*-algebra $A$ and a Hilbert
space $H$ completed with respect to the norm which is derived from the
$A$-valued inner product $\langle a \otimes h, b \otimes g \rangle = ab^*
\langle h,g \rangle_H$ for $a,b \in A$, $h,g \in H$. It is denoted by $A
\otimes H$. This kind of Hilbert C*-modules is very important: every
finitely generated Hilbert $A$-module $\mathcal M$ can be embedded into some
Hilbert $A$-module $A^n = A \odot {\mathbb C}^n$ for finite $n \in \mathbb N$
as an orthogonal summand, and every countably generated Hilbert $A$-module
can be realized as an orthogonal summand of the Hilbert $A$-module $l_2(A) =
A \otimes l_2$ in such a way that its orthogonal complement is isometrically
isomorphic to $A \otimes l_2$ again, (Kasparov's theorem\cite{Lance,NEWO}).
To name two further examples, the set of all continuous sections of a certain
vector bundle over a base space $X$ becomes naturally a finitely generated
Hilbert C($X$)-module, and conditional expectations on C*-algebras $B$ turn
these C*-algebras into pre-Hilbert $A$-modules over the image C*-subalgebra
$A \subseteq B$. (By the way, the interrelation between vector bundles and
finitely generated projective C*-modules over commutative C*-algebras is
a categorical equivalence similar to Gel'fand's theorem for C*-algebras, see
Serre-Swan's {theorem\cite{NEWO}}.)

The reader should be aware that the theory of Hilbert C*-modules has much
more exceptional examples in comparison to Hilbert space theory then one
can think of. E.g.~the analog of the Riesz' representation theorem for
bounded $A$-linear functionals on Hilbert $A$-modules is not valid, in
general. Therefore, bounded modular operators may be non-adjointable,
norm-closed Hilbert $A$-submodules can lack the property to be direct
summands, the notions of topological and orthogonal direct summands are
different, and more. There are even examples of non-countably generated
Banach C*-modules which can be turned into Hilbert C*-modules in at least two
ways, however the two C*-valued inner products are not unitarily isomorphic,
and so the resulting two norms give rise to non-isometrically isomorphic
Banach C*-modules, (\cite{Frank:99}). Luckily, the latter oddity cannot
appear for countably generated Hilbert C*-modules. For a concrete example
consider the C*-algebra $A={\rm C}([0,1])$of all continuous functions on the
unit interval, its ideal $I={\rm C}_0((0,1])$ of all continuous functions
vanishing at zero and the Hilbert $A$-module $\mathcal M = A \oplus I$ with
the standard $A$-valued inner product inherited from $\mathcal M \subset A^2$.
The operator $T: (a,i) \to (i,0)$ is non-adjointable, and the Hilbert
$A$-submodule $\mathcal N = \{ (i,i) : i \in I \}$ is a topological direct
summand, but not an orthogonal one.

In the light of these circumstances the results on the existence and on the
properties of modular frames for finitely or countably generated Hilbert
C*-modules presented below become the more remarkable. We start with a
definition of modular frames which takes an inequality in the positive cone
of the C*-algebra of coefficients as its initial point.

\begin{definition}
  Let $A$ be a unital C*-algebra and ${\mathbb J}$ be a finite or countable
  index set.
  A sequence $\{ x_j : j \in {\mathbb J} \}$ of elements in a Hilbert
  $A$-module $\mathcal M$ is said to be a {\it frame} if there are real
  constants $C,D > 0$ such that
  \begin{equation} \label{ineq-frame}
      C \cdot \langle x,x \rangle \leq
      \sum_{j} \langle x,x_j \rangle \langle x_j,x \rangle \leq
      D \cdot \langle x,x \rangle
  \end{equation}
  for every $x \in \mathcal M$. The optimal constants (i.e.~maximal for $C$
  and minimal for $D$) are called {\it frame bounds}. The frame $\{ x_j :
  j \in {\mathbb J} \}$ is said to be a {\it tight frame} if $C=D$, and said
  to be {\it normalized} if $C=D=1$.
  We consider {\it standard} (normalized tight) frames in the main for which
  the sum in the middle of the inequality (\ref{ineq-frame}) always converges
  in norm. For non-standard frames the sum in the middle converges only weakly
  for at least one element of $\mathcal M$.

  \smallskip
  A sequence $\{ x_j \}_j$ is said to be a {\it standard Riesz basis of
  $\mathcal M$} if it is a standard frame and a generating set with the
  additional property that $A$-linear combinations $\sum_{j \in S} a_jx_j $
  with coefficients $\{ a_j \} \in A$ and $S \in {\mathbb J}$ are equal to zero
  if and only if in particular every summand $a_jx_j$ equals zero for $j \in
  S$. A generating sequence $\{ x_j \}_j$ with the described additional
  property alone is called a {\it Hilbert basis of $\mathcal M$}.

  An {\it inner summand of a standard Riesz basis} of a Hilbert $A$-module
  $\mathcal L$ is a sequence $\{ x_j \}_j$ in a Hilbert $A$-module $\mathcal M$
  for which there is a second sequence $\{ y_j \}_j$ in a Hilbert $A$-module
  $\mathcal N$ such that ${\mathcal L} \cong {\mathcal M} \oplus {\mathcal N}$
  and the sequence consisting of the pairwise orthogonal sums $\{ x_j
  \oplus y_j \}_j$ in the Hilbert $A$-module ${\mathcal M} \oplus {\mathcal N}$
  is the initial standard Riesz basis of $\mathcal L$.

  \smallskip
  Two frames $\{ x_j \}_j$, $\{ y_j \}_j$ of Hilbert $A$-modules $H_1$, $H_2$,
  respectively, are {\it unitarily equivalent} (resp., {\it similar}) if there
  exists a unitary  (resp., invertible adjointable) linear operator $T: H_1 \to
  H_2$ such that $T(x_j)=y_j$ for every $j \in \mathbb J$.
\end{definition}

Analyzing this definition we do not know whether a frame is a generating set,
or not. This will turn out to hold only during our investigations.
We observe that for every (normalized tight) frame $\{ x_i \}_i$ of a
Hilbert space $H$ the sequence $\{ 1_A \otimes x_i \}_i$ is a standard
(normalized tight) module frame of the Hilbert $A$-module $\mathcal M =
A \otimes H$ with the same frame bounds. So standard modular frames exist in
abundance in the canonical Hilbert $A$-modules. At the same time wavelet
theorists see that the C*-algebra $A$ opens up an additional degree of freedom
for constructions and investigations. For the existence of
standard modular frames in arbitrary finitely or countably generated Hilbert
$A$-modules we obtained the following simple fact:

\begin{theorem} {\rm (\cite{FL:98,FL:00})} \newline 
  For every $A$-linear partial isometry $V$ on $A^n$ (or $l_2(A)$) the image
  sequence $\{ V(e_j) \}_j$ of the standard orthonormal basis $\{ e_j \}_j$
  is a standard normalized tight frame of the image $V(A^n)$ (or $V(l_2(A))$).
  Consequently, every algebraically finitely generated or countably generated
  Hilbert $A$-module $\mathcal M$ possesses a standard normalized tight frame
  since they can be embedded into these standard Hilbert $A$-modules as
  orthogonal summands.
\end{theorem}

\begin{problem}
  Does every Hilbert C*-module admit a modular frame ?
\end{problem}

The main property of frames for Hilbert spaces is the existence of the
reconstruction formula that allows a simple standard decomposition of every
element of the spaces with respect to the frame. We found that almost all
the results for the Hilbert space situation described in (\cite{HL}) can be
recovered. Sometimes the way of proving is exceptional long, for example to
show that modular Riesz bases $\{ x_i \}_i$ that are normalized tight frames
have to be orthogonal bases for which the values $\{ \langle x_i,x_i \rangle
\}_i$ are all projections. Let us first formulate the reconstruction formula
for normalized tight frames without the restriction to be standard:

\begin{theorem}  {\rm (Th.~4.1\cite{FL:98} )} \newline
  Let $A$ be a unital C*-algebra, $\mathcal M$ be a finitely or countably
  generated Hilbert $A$-module and $\{ x_j \}_j$ be a normalized tight frame
  of $\mathcal M$. Then the {\it reconstruction formula}
   \begin{equation}  \label{eq1}
     x= \sum_j \langle x,x_j \rangle x_j
   \end{equation}
  holds for every $x \in \mathcal M$ in the sense of convergence w.r.t.~the
  topology that is induced by the set of semi-norms 
  $\{ |f(\langle .,. \rangle )|^{1/2} \, : \, f \in A^* \}$.
  The sum converges always in norm if and only if the frame $\{ x_j \}_j$ is
  standard.      \newline
  Conversely, a finite set or a sequence $\{ x_j \}_j$ satisfying the formula
  (\ref{eq1}) for every $x \in \mathcal M$ is a normalized tight frame of
  $\mathcal M$.
\end{theorem}

For a proof we have to refer to (\cite{FL:98,FL:00}) since the proof is to
long to be reproduced here, and the statement for non-standard normalized
tight frames is some kind of summary of the entire work done. With the
experience on the possible oddities of Hilbert C*-module theory in comparison
to Hilbert space theory the following crucial fact is surprising because of
the generality in which it holds. The existence and the very good properties
of the frame transform of standard frames give the chance to get far reaching
results analogous to those in the Hilbert space situation. Again, the proof
is more complicated than the known one in the classical Hilbert space case,
cf.~(\cite{HL}).

\begin{theorem}   {\rm (Th.~4.2\cite{FL:98} )}   \newline
  Let $A$ be a unital C*-algebra, $\mathcal M$ be a finitely or countably
  generated Hilbert $A$-module and $\{ x_j \}_j$ be a standard frame
  of $\mathcal M$. The {\it frame transform of the frame $\{ x_j \}_j$} is
  defined to be the map
    \[
       \theta: \mathcal M \rightarrow l_2(A) \quad , \qquad \theta(x)
       = \{ \langle x,x_j \rangle \}_j 
    \]
  that is bounded, $A$-linear, adjointable and fulfills $\theta^*(e_j)=x_j$
  for a standard orthonormal basis $\{ e_j \}_j$ of the Hilbert $A$-module
  $l_2(A)$ and all $j \in \mathbb J$.   \newline
  Moreover, the image $\theta(\mathcal M)$ is an orthogonal summand of
  $l_2(A)$. For normalized tight frames we additionally get $P(e_j)=\theta
  (x_j)$ for any $j \in \mathbb J$, and $\theta$ is an isometry in that case.
\end{theorem}

The frame transform $\theta$ is the proper tool for the description of
standard frames.

\begin{theorem}  {\rm (Th.~6.1\cite{FL:98} )} \newline
  Let $A$ be a unital C*-algebra, $\mathcal M$ be a finitely or countably
  generated Hilbert $A$-module and $\{ x_j \}_j$ be a standard frame
  of $\mathcal M$. Then the {\it reconstruction formula}
    \[
         x= \sum_j \langle x,S(x_j) \rangle x_j
    \]
  holds for every $x \in \mathcal M$ in the sense of norm convergence, where
  the operator $S:=(\theta^*\theta)^{-1}$ is positive and invertible.
\end{theorem}

The sequence $\{ S(x_j) \}_j$ is a standard frame again, the {\it canonical
dual frame of} $\{ x_j \}_j$. The operator $S$ is called the {\it frame
operator of $\{ x_j \}_j$} on $\mathcal M$. The two theorems on reconstruction
above bring to light some key properties of frame sequences:

\begin{corollary} {\rm (Th.~5.4, 5.5\cite{FL:98}) }  \newline
  Every standard frame of a finitely or countably generated Hilbert $A$-module
  is a set of generators.  \newline
  Every finite set of algebraic generators of a finitely generated Hilbert
  $A$-module is a (standard) frame.  
\end{corollary}

\begin{corollary}
  Every standard frame can be realized as the image of a standard normalized
  tight frame by an adjointable invertible bounded module operator. In
  particular, every standard Riesz basis is the image of an orthogonal Hilbert
  basis $\{ x_i \}_i$ with projection-valued values $\{ \langle x_i,x_i \rangle
  \}_i$ by an adjointable invertible bounded module operator.
\end{corollary}

Let us show that there exists countably generated Hilbert C*-modules without
standard Riesz bases: let $A={\rm C}([0,1])$ be the C*-algebra of all
continuous functions on the unit interval $[0,1]$ and consider its ideal and
Hilbert $A$-submodule $\mathcal M = {\rm C}_0([0,1])$ of all continuous
functions vanishing at zero. The module $\mathcal M$ is countably generated
by the Stone-Weierstrass theorem. So it admits normalized tight frames.
However, if it would have a standard Riesz basis, then it would contain an
orthogonal Hilbert basis $\{ x_i \}_i$ with non-trivial projection-valued
values $\{ \langle x_i,x_i \rangle \}_i$. At the other side, the only
projection contained in the range of the $A$-valued inner product is the zero
element. So it does not contain any standard Riesz basis.

We close our considerations on the frame transform arising from standard
frames with a statement on the relation between unitary equivalence (or
similarity) of frames and the characteristics of the image of the frame
transform $\theta$. 

\begin{theorem} {\rm (\cite[Th.~7.2]{FL:98}) } \newline
  Let $A$ be a unital C*-algebra and $\{ x_j \}_j$ and $\{ y_j \}_j$ be
  standard (normalized tight) frames of Hilbert $A$-modules ${\mathcal M}_1$
  and ${\mathcal M}_2$, respectively. The following conditions are
  equivalent:
  \begin{list}{(\roman{marke})}{\usecounter{marke}}
   \item The frames $\{ x_j \}_j$ and $\{ y_j \}_j$ are unitarily equivalent
         or similar.
   \item Their frame transforms $\theta_1$ and $\theta_2$ have the same range
         in $l_2(A)$.
   \item The sums $\sum_j a_j x_j$ and $\sum_j a_j y_j$ equal zero for exactly
         the same Banach A-submodule of sequences $\{ a_j \}_j$ of $l_2(A)$.
  \end{list}
\end{theorem}

These results look pretty much the same as for the Hilbert space situation
what makes them easy to apply. However, the proofs are more difficult, and
they required an extensive search for possible counterexamples and obstacles
beforehand. There are also statements that do not transfer to Hilbert
C*-modules. For example, standard Riesz bases of Hilbert C*-modules may have
more than one dual frame because of the existence of zero-divisors in the
C*-algebra of coefficients, see Example {6.4\cite{FL:98}} and Corollary
{6.6\cite{CHL,HL}}.

More results on disjointness and inner sums of frames, as well as on various
kinds of frame decompositions can be found in (\cite{FL:00}). Our basic
publications\cite{FL:98,FL:00} also contain a number of illustrating examples
and counterexamples we would like to refer to.

%%%%%%%%%%%%%%%%%%%%%%%%%%%%%%%%%%%%%%%%%%%%%%%%%%%%%%%%%%%%%%%%%%%%%%%%%%%%%

\section{INVARIANTS OF FINITELY GENERATED PROJECTIVE C*-MODULES}

\noindent
The best way for the description of the inner structure of Hilbert spaces is
the selection of a ((ortho-) normal) basis and the characterization of elements
by their coordinates. The notion of a basis makes essentially use of the
notion of linear independence of vectors. Turning to finitely generated
projective C*-modules over a certain fixed unital C*-algebra $A$ we are most
often faced with the absence of a reasonable notion of '$A$-linear independence'
of sets of module elements, e.g.~finite sets of algebraic generators. Also,
very often we have a lot of non-isomorphic $A$-modules possessing generating
sets with the same number of algebraic generators. The task is to find
additional invariants for the distinction of projective $A$-modules in such
situations.

Fortunately, any set of algebraic generators of a finitely generated projective
C*-module $\mathcal M$ over a unital C*-algebra $A$ is a frame, a fact
shown by the authors\cite{FL:98} in 1998. Furthermore, for every set of algebraic
generators $\{ x_1, ..., x_k \}$ of $\mathcal M$ there exists an $A$-valued
inner product $\langle .,. \rangle$ on $\mathcal M$ turning this set 
$\{ x_1, ..., x_k \}$ into a normalized tight frame. We show that the
knowledge of the values $\{ \langle x_i,x_j \} : 1 \leq i \leq j \leq k \}$
turns out to be sufficient to describe the $A$-module $\mathcal M$ up to
uniqueness. Note that the elements $\{ x_1,...,x_k \}$ need not to be
($A$-)linearly independent, in general.

\begin{theorem}
  Let $A$ be a unital C*-algebra and let $\{ \mathcal M, \langle .,.
  \rangle_{\mathcal M} \}$ and $\{ \mathcal N, \langle .,. \rangle_{\mathcal N}
  \}$ be two finitely generated Hilbert $A$-modules. Then the following
  conditions are equivalent:
   \begin{list}{(\roman{marke})}{\usecounter{marke}}
    \item  $\{ \mathcal M, \|.\|_{\mathcal M} \}$ and $\{ \mathcal N,
       \|.\|_{\mathcal N} \}$ are isometrically isomorphic as Banach
       $A$-modules.
    \item $\{ \mathcal M, \langle .,. \rangle_{\mathcal M}$ and $\{ \mathcal
       N, \langle .,. \rangle_{\mathcal N} \}$ are unitarily isomorphic as
       Hilbert $A$-modules.
    \item There are finite normalized tight frames $\{ x_1,...,x_k
       \}$ and $\{ y_1,...,y_l \}$ of $\mathcal M$ and $\mathcal N$,
       respectively, such that $k=l$, $x_i \not= 0$ and $y_i \not= 0$ for any
       $i=1,...,k$, and $\langle x_i,x_j \rangle_{\mathcal M} = \langle y_i,y_j
       \rangle_{\mathcal N}$ for any $1 \leq i \leq j \leq k$.
   \end{list}
\end{theorem}

\begin{proof} The equivalence of the conditions (i) and (ii) has been
shown for countably generated Hilbert $A$-modules in Theorem 4.1\cite{Frank:99}.
The implication (ii)$\to$(iii) can be seen to hold setting $y_i = U(x_i)$
for the existing unitary operator $U: \mathcal M \to \mathcal N$ and for
$i=1,...,k$. The demonstration of the inverse implication requires slightly
more work. For the given normalized tight frames $\{ x_1,...,x_k \}$
and $\{ y_1,...,y_l \}$ of $\mathcal M$ and $\mathcal N$, respectively, we
define a $A$-linear operator $V$ by the rule $V(x_i)=y_i$, $i=1,...,k$.
For this operator $V$ we obtain the equalities
  \begin{eqnarray*}
     \langle V(x),y_j \rangle_{\mathcal N} & = &
        \sum_{i=1}^k \langle V(x),y_i \rangle_{\mathcal N} \langle y_i,y_j
        \rangle_{\mathcal N} 
     = \sum_{i=1}^k \left\langle \sum_{m=1}^k \langle x,x_m \rangle_{\mathcal M}
        V(x_m),y_i \right\rangle_{\mathcal N} \langle x_i,x_j \rangle_{\mathcal M} \\
     &=& \sum_{i=1}^k \sum_{m=1}^k \langle x,x_m \rangle_{\mathcal M} \langle x_m,
        x_i \rangle_{\mathcal M} \langle x_i,x_j \rangle_{\mathcal M} 
     = \sum_{m=1}^k  \langle x,x_m \rangle_{\mathcal M} \left\langle x_m,
        \sum_{i=1}^k \langle x_j,x_i \rangle_{\mathcal M} x_i
        \right\rangle_{\mathcal M} \\
     &=& \left\langle x , \sum_{m=1}^k \langle x_j,x_m \rangle_{\mathcal M} x_m
        \right\rangle_{\mathcal M} 
     = \langle x,x_j \rangle_{\mathcal M}\\
  \end{eqnarray*}
which hold for every $x \in \mathcal M$. Consequently,
\[
  \langle x,x \rangle_{\mathcal M} =
  \sum_{i=1}^k \langle x,x_i \rangle_{\mathcal M} \langle x_i,x \rangle_{\mathcal M}=
  \sum_{i=1}^k \langle V(x),y_i \rangle_{\mathcal N} \langle y_i,V(x) \rangle_{\mathcal N}=
  \langle V(x),V(x) \rangle_{\mathcal N}
\]
for any $x \in \mathcal M$, and the operator $V$ is unitary.
\end{proof}

\begin{corollary}
Every finitely generated projective $A$-module $\mathcal M$ over a unital
C*-algebra $A$ can be reconstructed up to isomorphism from the following data:
   \begin{list}{(\roman{marke})}{\usecounter{marke}}
    \item   A finite set of algebraic non-zero modular generators $\{
        x_1,...,x_k \}$ of $\mathcal M$.
    \item   A symmetric $k \times k$ matrix $( a_{ij} )$ of elements from $A$,
        where $a_{ij}$ is supposed to be equal to $\langle x_i,x_j \rangle_0$
        for $1 \leq i \leq j \leq k$ and for the (existing and unique)
        $A$-valued inner product $\langle .,. \rangle_0$ on $\mathcal M$ that
        turns the set of algebraic modular generators $\{ x_1,...,x_k\}$ into
        a normalized tight frame of the Hilbert $A$-module $\{ \mathcal M,
        \langle .,. \rangle_0 \}$.
   \end{list}
The number of elements in sets of algebraic modular generators of $\mathcal M$
has a minimum, and it suffices to consider sets of generators of minimal
length. Then the modular invariants can be easier compared permuting the
elements of the generating sets if necessary.
\end{corollary}

\begin{proof}
We have already pointed out that the set of algebraic generators $\{ x_1,...
,x_k \}$ of $\mathcal M$ is a frame with respect to any $A$-valued inner
product on $\mathcal M$ which turns $\mathcal M$ into a Hilbert $A$-module.
That means the inequality
\[
  C \cdot  \langle x,x \rangle \leq \sum_{i=1}^k \langle x,x_i \rangle
  \langle x_i,x \rangle \leq D \cdot \langle x,x \rangle
\]
is satisfied for two finite positive real constants $C,D$ and any $x \in
\mathcal M$, see {Theorem 5.9\cite{FL:98}}. What is more, for any frame of
$\mathcal M$ there exists another $A$-valued inner product $\langle .,.
\rangle_0$ on $\mathcal M$ with respect to which it becomes normalized tight.
The latter inner product is unique as shown by Corollary 4.3,
Theorem 6.1\cite{FL:98} and {Theorem 4.4\cite{FL:00}}. So assertion (iv) of the
previous theorem demonstrates the complete assertion.
\end{proof}

\smallskip
We can say more in case the finitely generated Hilbert $A$-module contains
a modular Riesz basis, i.e.~a finite set of modular generators $\{ x_1,...,x_k
\}$ such that the equality $0= a_1x_1 + ... + a_kx_k$ holds for certain
coefficients $\{ a_1,...,a_k \} \subset A$ if and only if $a_ix_i=0$ for any
$i=1,...,k$. Obviously, a modular Riesz basis is minimal as a set of modular
generators, i.e.~we cannot drop any of its elements preserving the generating
property. However, there can exist totally different Riesz bases for the same
module that consist of less elements, cf.~{Example 1.1\cite{FL:98}}. Note
that the coefficients $\{ a_1,...,a_k \}$ can be non-trivial even if $a_ix_i=0$
for any index $i$ since every non-trivial C*-algebra $A$ contains zero-divisors.
Not every Hilbert C*-module with a normalized tight modular frame does
possess a modular Riesz basis. For an example we refer to {Example
2.4\cite{FL:00}}.

In case of finitely generated projective W*-modules (and therefore, in the case
of Hilbertian modules over finite W*-algebras) we are in the pleasant situation
that they always contain a modular Riesz basis by W.~L.~Paschke's Theorem
{3.1\cite{Pa}}. Moreover, by spectral decomposition every element $x$ of a
Hilbert W*-module $\mathcal M$ has a carrier projection of $\langle x,x
\rangle$ contained in the W*-algebra of coefficients $A$. So we can ascertain
the following fact:

\begin{proposition}
  Let $\mathcal M$ be a finitely generated projective $A$-module over a
  W*-algebra $A$ that possesses two finite modular Riesz bases $\{ x_1,...,x_k
  \}$ and $\{y_1,...,y_l \}$.
  Then there exists an $l \times k$ matrix ${\mathbf F} = (f_{ij})$, $i=1,...,l$,
  $,j=1,...,k$, with entries from $A$ such that $y_i = \sum_{j=1}^k f_{ij} x_j$
  for any $i = 1,...,l$, and analogously, there exists a $k \times l$ matrix
  ${\mathbf G} = (g_{ji})$ with entries from $A$ such that $x_j = \sum_{i=1}^l g_{ji} y_i$
  for any $j=1,...,k$.

  Suppose the left carrier projections of $f_{ij}$ and $g_{ji}$ equal the
  carrier projections of $\langle y_i,y_i \rangle$ and $\langle x_j,x_j \rangle$,
  respectively, and the right carrier projection of $f_{ij}$ and $g_{ji}$ equal
  the carrier projections of $\langle x_j,x_j \rangle$ and $\langle y_i,y_i
  \rangle$, respectively. Then the matrices $\mathbf F$ and $\mathbf G$ are
  Moore-Penrose invertible in $M_{kl}(A)$ and $M_{lk}(A)$, respectively. The
  matrix $\mathbf F$ is the Moore-Penrose inverse of $\mathbf G$, and vice
  versa.
\end{proposition}

\begin{proof}
Since both the modular Riesz bases are sets of modular generators of $\mathcal M$
we obtain two $A$-valued (rectangular, w.l.o.g.) matrices ${\mathbf F} =
(f_{ij})$ and ${\mathbf G} = (g_{ji})$ with $i=1,...,l$ and $j=1,...,k$ such
that
\[
  y_i = \sum_{m=1}^k f_{im}x_m \quad , \quad  x_j = \sum_{n=1}^l g_{jn} y_n
   \, .
\]
Combining these two sets of equalities in both the possible ways we obtain
\[
  y_i = \sum_{n=1}^l \left( \sum_{m=1}^k f_{im} g_{mn} \right) y_n  \quad , \quad
  x_j = \sum_{m=1}^k \left( \sum_{n=1}^l g_{jn} f_{nm} \right) x_m
\]
for $i=1,...,l$, $j=1,...,k$. Now, since we deal with sets of coefficients
$\{ f_{ij} \}$ and $\{ g_{ji} \}$ that are supposed to admit special carrier
projections,  the coefficients in front of the elements $\{ y_n \}$ and $\{
x_m \}$ at the right side can only take very specific values:
\[
   \sum_{m=1}^k f_{im} g_{mn} = \delta_{in} \cdot q_n \quad , \quad
   \sum_{n=1}^l g_{jn} f_{nm} = \delta_{jm} \cdot p_m  \, ,
\]
where $\delta_{ij}$ is the Kronecker symbol, $p_m \in A$ is the carrier
projection of $\langle x_m,x_m \rangle$ and $q_n \in A$ is the carrier
projection of $\langle y_n,y_n \rangle$. So ${\mathbf F} \cdot \mathbf G$ and
${\mathbf G} \cdot \mathbf F$ are positive idempotent diagonal matrices with
entries from $A$. The Moore-Penrose relations ${\mathbf F} \cdot {\mathbf G}
\cdot {\mathbf F} = \mathbf F$, ${\mathbf G} \cdot {\mathbf F} \cdot {\mathbf G}
= \mathbf G$, $({\mathbf F} \cdot {\mathbf G})^* = {\mathbf F} \cdot \mathbf G$
and $({\mathbf G} \cdot {\mathbf F})^* = {\mathbf G} \cdot \mathbf F$ turn out
to be fulfilled.
\end{proof}

In total we found a convinient way to characterize finitely generated
C*-modules over unital C*-algebras up to modular isomorphism by a small amount
of additional elements of the C*-algebra of coefficients derived from the
set of algebraic generators and from the module structure.

%%%%%%%%%%%%%%%%%%%%%%%%%%%%%%%%%%%%%%%%%%%%%%%%%%%%%%%%%%%%%%%%%%%%%%%%%%%%%

\section{AN OPERATOR-THEORETIC PROBLEM RESOLVED USING FRAME THEORY}

\noindent
One of the classical problems of operator theory is the following:
given a (finite or infinite) sequence $\{ b_i \}_{i=1}^\infty$ of bounded
operators on a certain separable Hilbert space $H$ fulfilling the equality
${\rm id}_{H} = \sum_{i=1}^\infty b_i^*b_i$, determine the nature of the operators
$\{ b_i \}_{i=1}^\infty$. A first account was found by R.~V.~Kadison
and J.~R.~Ringrose in II.11.2.24\cite{KaRi} using dilation methods,
i.e.~enlarging the Hilbert space $H$ and subsequently extending the operators.
We will use modular frame methods to resolve this problem without changing
the Hilbert space $H$. That way we demonstrate the power of the developed
methods. Some particular examples complete the picture.

\begin{proposition}
   Let $\{ b_i \}_{i=1}^\infty \in B(l_2)$ be a sequence of elements with
   the property that ${\rm id}_{l_2} = \sum_{i=1}^\infty b_i^*b_i$ in the sense
   of weak convergence. (In particular, the sequence could be norm-convergent,
   or only finitely many elements could be unequal to zero.)
   Then there exists a projection $p \in B(l_2)$, $p \sim {\rm id}_{l_2}$ via
   $uu^*=p$, $u^*u={\rm id}_{l_2}$, and a sequence of partial isometries
   $\{ v_i \}_{i=1}^\infty \in B(l_2)$ such that:
   \begin{list}{(\roman{marke})}{\usecounter{marke}}
    \item  $v_iv_i^*={\rm id}_{l_2}$, $v_i^*v_i$ are projections in $B(l_2)$ similar
       to ${\rm id}_{l_2}$, $(v_i^*v_i) \bot (v_k^*v_k)$ for any $i \not= k$,
       $\sum_{i=1}^\infty v_i^*v_i = {\rm id}_{l_2}$;
    \item the equality $b_i = u_i (v_iu)$ is valid for every $i \in
       \mathbb N$ and partial isometries $\{ u_i \}_{i=1}^\infty \subset B(l_2)$,
       each $u_i$ connecting the left carrier projections of $v_iu$ and of
       $b_i$, respectively.
   \end{list}
\end{proposition}

\begin{proof}
Let $A=B(l_2)$ be the set of all linear bounded operators on the Hilbert
space $l_2$, and recall that $l_2 = \oplus_{i=1}^\infty (l_2)_{(i)}$ as a
direct sum of copies of the Hilbert space $l_2$ itself. Consequently, there
are projections $\{ p_i \}_{i=1}^\infty$ such that $p_i \bot p_j$ for any
$i \not= j$, $\sum_{i=1}^\infty p_i = {\rm id}_{l_2}$, and $p_i \sim {\rm id}_{l_2}$ via
partial isometries $\{ v_i\}_{i=1}^\infty$ with $v_iv_i^* = {\rm id}_{l_2}$ and
$v_i^*v_i = p_i$.

Define $l_2(A)' := \{ \{a_i \}_{i=1}^\infty : \sup_{N} \left\| \sum_{i=1}^N
a_ia_i^* \right\| < \infty \}$, i.e. the set of all sequences of elements of
$A$ for which the series converges weakly in $A=B(l_2)$. It can be identified
with the set of all bounded $A$-linear maps $r$ of the Hilbert $A$-module
$l_2(A)$ into $A$ setting $a_i =r(e_i)$, where $e_i = (0,...,0,1_{A,(i)}, 0,...)$.
Since $A=B(l_2)$ is a W*-algebra, the set $l_2(A)'$ becomes a self-dual
Hilbert $A$-module ({\cite{Pa}}).

In fact, $l_2(A)'$ is isometrically isomorphic to $A=B(l_2)$ as a Hilbert
$A$-module. To see this fix an orthonormal basis $\{ e_i \}_{i=1}^\infty$ of
$l_2(A)'$ and an orthonormal basis $\{ v_i \}_{i=1}^\infty$ of $A$ and define
\[
   l_2(A)' \to A \quad, \qquad \{ a_i\}_{i=1}^\infty \to \sum_{i=1}^\infty
   a_i v_i \, ;
\]
\[
   A \to l_2(A)' \quad , \qquad a \to \{ av_i^* \}_{i=1}^\infty \, .
\]
The freedom of choice for this isomorphism is the careful selection of the
orthonormal bases of the self-dual Hilbert $A$-modules $l_2(A)'$ and $A$.

If the sequence $\{ b_i \}_{i=1}^\infty \in B(l_2)$ is given as described
above then  it is a (possibly non-standard) normalized tight frame of the
self-dual Hilbert $A$-module $A=B(l_2)$ since we have
\[
   a = a \cdot {\rm id}_{l_2} = a \sum_{i=1}^\infty  b_i^*b_i
   = \sum_{i=1}^\infty \langle a, b_i \rangle_{B(l_2)} b_i
\]
for every $a \in A = B(l_2)$. There exists a frame transform $\theta: B(l_2)
\to l_2(B(l_2))'$ defined by $a \to \{ \langle a,b_i \rangle_{B(l_2)}
\}_{i=1}^\infty$. The image of $\theta$ is a direct summand and self-dual
Hilbert $A$-submodule of $l_2(A)'$. Moreover, $\theta$ is an isometry.
Continuing the isometry $\theta$ to an isometry $\theta': B(l_2) \to B(l_2)$
using the isometric isomorphism between $l_2(A)'$ and $A$ we obtain
\[
   \theta' : B(l_2) \to B(l_2) \quad , \qquad a \to \sum_{i=1}^\infty \langle
   a,b_i \rangle_{B(l_2)} v_i = \sum_{i=1}^\infty ab_i^*v_i
\]
for any $a \in A$. The structure of direct orthogonal summands of the
Hilbert $B(l_2)$-module $B(l_2)$ is well-known, they are all generated by
multiplying $B(l_2)$ by a specific orthogonal projection from the right.
So we can characterize the image of $\theta'$ in $B(l_2)$ as $B(l_2) p$ for
some $p=p^2 \geq 0$ of $B(l_2)$. Note that $\theta'({\rm id}_{l_2})=p$ since the
module generator ${\rm id}_{l_2}$ is mapped to the module generator $p$.
Let $u \in B(l_2)$ be the isometry linking $p$ to ${\rm id}_{l_2}$ with $uu^*=p$,
$u^*u={\rm id}_{l_2}$. Now, we establish some information on the adjoint operator
of $\theta'$:
\begin{eqnarray*}
   \langle a,b_i \rangle_{B(l_2)}
   & = & \langle a, (\theta')^*(v_i) \rangle_{B(l_2)} 
     =   \langle \theta'(a), v_i \rangle_{B(l_2)} \\
   & = & \left( \sum_{j=1}^\infty \langle a,b_j \rangle_{B(l_2)}
         v_j \right) v_i^* 
     =   \left( \sum_{j=1}^\infty \langle a,b_j \rangle_{B(l_2)} v_j \right)p
         v_i^* \\
   & = & \left\langle  a , \sum_{j=1}^\infty v_i p v_j^*b_j \right\rangle_{
         B(l_2)} \, .
\end{eqnarray*}
Consequently, we get the frame decomposition $b_i = \sum_{j=1}^\infty v_i p
v_j^*b_j $ for any $i \in \mathbb N$. The frame coefficients $\{ v_i p v_j^*
\}_{j=1}^\infty$ may be not the optimal ones. By Prop.~6.6\cite{FL:98} we
have the general inequality:
\begin{equation} \label{eq-1}
   (v_iu)(v_iu)^* =
   v_ipv_i^* = \sum_{j=1}^\infty v_ipv_j^*v_jpv_i^* \geq  \sum_{j=1}^\infty
   \langle b_i,b_j \rangle \langle b_j,b_i \rangle = b_ib_i^*
\end{equation}
for any $i \in \mathbb N$. By Theorem 2.1\cite{Ped} we obtain
$b_i = u_i (v_iu)$ for some elements $\{ u_i \}_{i=1}^\infty$ with
$\| u_i \| \leq 1$ and
\[
   u_i = {\rm strong-lim} \,\, b_i^* (\varepsilon + v_ipv_i^*)^{-1}v_iu =
   {\rm strong-lim} \,\,  b_i^* v_iu(\varepsilon + u^*p_iu)^{-1} \, .
\]
Of course the root of $v_ipv_i^*$ seems to be selected in an artificial way,
the element $v_ip$ would do the job as well. However, the following inequality
gives some more information on the background of the choice made:
\[
   {\rm id}_{l_2} = \sum_{i=1}^\infty b_i^*b_i
        =     \sum_{i=1}^\infty u^*v_i^*u_i^*u_iv_iu 
        \leq  \sum_{i=1}^\infty u^*v_i^* \|u_i\|^2 v_iu 
        \leq  u^* \left( \sum_{i=1}^\infty p v_i^*v_i p \right) u
        =     {\rm id}_{l_2} \, .
\]
Since the left end equals the right end and $b_i^*b_i \leq (v_iu)^*(uv_i)$
holds for every $i \in \mathbb N$ the equality $b_i^*b_i = (v_iu)^*(v_iu)$
turns out to be valid for every $i \in \mathbb N$.  Consequently, the linking
elements $\{ u_i \}_{i=1}^\infty$ can be selected as partial isometries of
$B(l_2)$ mapping the left carrier projection of $v_iu$ to the left carrier
projection of $b_i$ for each $i \in \mathbb N$ because $B(l_2)$ is a von
Neumann algebra and any root of a given positive operator can be described
this way. By the inequality (\ref{eq-1}) the left carrier projection of $b_i$
has to be lower-equal than the left carrier projection of $v_iu$ for any
$i \in \mathbb N$.
\end{proof}

\begin{example}        {\rm
  Suppose, the set of operators $\{ b_i \}_{i=1}^\infty$ is a set of pairwise
  orthogonal (positive) projections $\{ p_i \}_{i=1}^\infty$ defined on $l_2$.
  This is the simplest situation one can think of.
  Then $u=p={\rm id}_{l_2}$ and $v_i=u_i =p_i$ for any index $i$.  \newline
  Furthermore, there is the classical situation of generators of Cuntz algebras
  ${\mathcal O}_n$ and ${\mathcal O}_\infty$, where the operators
  $\{ b_i \}_{i=1}^\infty$ are partial isometries themselves, with the
  properties required above.
  Here $u=p={\rm id}_{l_2}$ and $v_i=b_i$ for any index $i$. The partial isometries
  $u_i$ equal to the left carrier projections of the partial isometries $b_i$
  in the given situation.  }
\end{example}

Generally speaking, the crucial rule is played by the projection $p$
corresponding to the sequence $\{ b_i \}_{i=1}^\infty$ via its frame transform,
and by its partition $\{ p_ip \}_{i=1}^\infty$ by a chain of pairwise orthogonal
projections $\{ p_i \}_{i=1}^\infty$ summing up to one and each being similar
to the identity operator on $l_2$.

%%%%%%%%%%%%%%%%%%%%%%%%%%%%%%%%%%%%%%%%%%%%%%%%%%%%%%%%%%%%%%%%%%%%%%%%%%%%%

\section{APPROXIMATION OF FRAMES BY (NORMALIZED) TIGHT ONES}

\noindent
In the present section we consider a question on the approximation of frames
of Hilbert spaces $H$ by (normalized) tight ones that is related to certain
methods of orthogonalization and renormalization of Hilbert bases. We have to
resort to Hilbert spaces instead of Hilbert C*-modules since the problem is
too complex to be treated in full generality. We give comments on the more
general setting whereever possible. Unfortunately, we did not find a final
solution of the problem, rather we obtained hints to the complexity and
difficulty of it. The partial results are nevertheless worth to be discussed.
Generally speaking, most distance measures on sets of frame operators seem to
have rather an $L^\infty$-character than an $L^2$-character, exept the
Hilbert-Schmidt norm. So the stressed for uniqueness of best approximating
tight frames often cannot be obtained.

The following fundamental problem has been pointed out by R.~Balan and the
first author\cite{BF} in July 1999 summarizing earlier investigations:

\begin{problem}     {\rm
  Are there distance measures on the set of frames of Hilbert subspaces $K$
  of $H$ with respect to which a multiple of the normalized tight frame
  $\{ S^{1/2} (x_i) \}_i$ is the closest (normalized) tight frame to the
  given frame $\{ x_i \}_i$ of the Hilbert subspace $K \subseteq H$, or at
  least one of the closest (normalized) tight frames?
\newline 
  If there are other closest (normalized) tight frames with respect to the
  selected distance measures, do they span the same Hilbert subspaces of $H$?
  If not, how are the positions of these subspaces with respect to $K \subseteq
  H$?    }
\end{problem}

Let $\{ e_i \}_i$ be an orthonormal basis of the Hilbert space $l_2(\mathbb C)$.
For the analysis operator $T: K \to l_2(\mathbb C)$, $T(x_i) = \{ \langle x_i,
e_i \rangle \}_i$, of a subspace frame $\{ x_i \}_i$ of $K$ there is a polar
decomposition $T = V S^{-1/2}$, where $S=(T^*T)^{-1}$ denotes the frame
operator. We can easily check that $S^{1/2}(x_i) = V^*(e_i)$ for any $i \in
\mathbb N$. Also the projection $P: l_2(\mathbb C) \to TT^*(l_2(\mathbb C))$
plays an important role.

\smallskip
Looking through the literature there are two approaches to this
problem, one due to R.~Balan\cite{Balan} and the other due to V.~I.~Paulsen,
T.~R.~Tiballi\cite{Tib} and the first author\cite{FPT}.
R.~Balan\cite{Balan} starts with the definition of three distance measures
for pairs of frames:
The frame $\{ x_i \}_i$ of the Hilbert space $H$ is said to be {\it
quadratically close} to the frame $\{ y_i \}_i$ of $H$ if there exists a
non-negative number $C$ such that the inequality
\[
   \left\| {\sum}_i c_i (x_i-y_i) \right\|
   \leq C \cdot \left\| {\sum}_i c_iy_i \right\|
\]
is satisfied.  The infimum of all such constants $C$ is denoted by $c(y,x)$.
In general, if $C \geq c(y,x)$ then $C(1-C)^{-1} \geq c(x,y)$, however this
distance measure is not reflexive. Two frames $\{ x_i \}_i$ and $\{ y_i \}_i$
of a Hilbert space $H$ are said to be {\it near} if  $d(x,y)= \log (\max
(c(x,y),c(y,x))+1) < \infty$. They are near if and only if they are similar,
(Th.~2.4\cite{Balan}). The distance measure $d(x,y)$ is an equivalence
relation and fulfils the triangle inequality.

\begin{theorem} {\rm (Th.~3.1\cite{Balan})} \newline
  For a given frame $\{ x_i \}_i$ of $H$ with frame bounds $C,D$ the distance
  measures admit their infima at
  \[
   \min \, c(y,x) = \min \, c(x,y) = \frac{\sqrt{D}-\sqrt{C}}{\sqrt{D}+\sqrt{C}}
   \, , \,\, \min \, d(x,y) = \frac{1}{4} (\log (D) - \log (C)) \, .
  \]
  These values are achieved by the tight frames
  \begin{equation} \label{eqn1}
   \left\{ \frac{\sqrt{C}+\sqrt{D}}{2} S^{1/2}(x_i) \right\}_i \, , \,
   \left\{ \frac{2\sqrt{CD}}{\sqrt{C}+\sqrt{D}} S^{1/2}(x_i) \right\}_i \, ,
   \{ \sqrt[4]{CD} S^{1/2}(x_i) \}_i   \, ,
  \end{equation}
  listed in the same order as the measures above. The solution may not be
  unique, in general, however any tight frame $\{ y_i \}_i$ of $H$ that
  achieves the minimum of one of the three distance measures $c(y,x)$, $c(x,y)$
  and $d(x,y)$ is unitarily equivalent to the corresponding solutions listed
  above.
\end{theorem}

The difference of the connecting unitary operator and the product of minimal
distance times either $S^{1/2}$ or $S^{-1/2}$ fulfils a certain measure-specific
operator norm equality which can be found at {Th.~3.1\cite{Balan}}. We 
point out that the first constant at (\ref{eqn1}) is the arithmetic mean
of $\sqrt{C}$ and $\sqrt{D}$, the second one is their harmonic mean and the
third one is their geometrical mean.

The results by V.~I.~Paulsen, T.~R.~Tiballi\cite{Tib} and the first author\cite{FPT}
are of slightly different character, however the operator $(P-|T^*|)$ has
to be Hilbert-Schmidt for their validity.

\begin{theorem} {\rm (Th.~2.3\cite{FPT})} \newline
  The operator $(P-|T^*|)$ is Hilbert-Schmidt if and only if the sum
  ${\sum}_{j=1}^\infty \| \mu_j-x_j \|^2$ is finite for at least one normalized
  tight frame $\{ \mu_i \}_i$ of a Hilbert subspace $L$ of $H$
  that is similar to $\{ x_i \}_i$. In this situation the estimate
  \[
     {\sum}_{j=1}^\infty \| \mu_j-x_j \|^2 \geq
     {\sum}_{j=1}^\infty \| S^{1/2}(x_j)-x_j \|^2 = \|(P-|T^*|)\|^2_{c_2}
  \]
  is valid for every normalized tight frame $\{ \mu_i \}_i$
  of any Hilbert subspace $L$ of $H$ that is similar to $\{ x_i \}_i$,
  where $\| \cdot \|_{c_2}$ denotes the Hilbert-Schmidt norm.
  (The left sum can be infinite for some choices of  subspaces $L$ and
  normalized tight frames $\{ \mu_i \}_i$ for them.)
\newline
  Equality appears if and only if $\mu_i =S^{1/2}(x_i)$ for any $i \in \mathbb
  N$. Consequently, {\it the symmetric approximation of a frame $\{ x_i \}_i$
  in a Hilbert space $K \subseteq H$} is the normalized tight frame
  $\{ S^{1/2}(x_i) \}_i$ spanning the same Hilbert subspace $L \equiv K$ of $H$
  and being similar to $\{ x_i \}_i$ via the invertible operator $S^{-1/2}$.
\end{theorem}

Applying the theorem to appropriate Riesz bases $\{ x_i \}_i$ the normalized
tight frame $\{ S^{1/2}(x_i) \}_i$ turns out to be the symmetric
orthogonalization of this basis as discovered by P.-O.~L\"owdin\cite{Loew} in
1948. This is why the denotation `symmetric approximation' has been selected
for the normalized tight frame $\{ S^{1/2}(x_i) \}_i$.

The resultsin ({\cite{FPT}}) generalize to modular frames of countably generated
Hilbert C*-modules over commutative C*-algebras, whereas for respective Hilbert
C*-modules over non-commutative C*-algebras the proofs cannot reproduced at
several crucial places where commutativity of the C*-algebra of coefficients
is essential. So the formulation of propositions on the non-commutative case
is an open problem at present.

\smallskip
To overcome the difficulties with the non-commutativity of the C*-algebra of
coefficients we consider distance-measures based on the various frame operators.
The properties of these operators do not depend on the choice of the set of
coefficients ({\cite{FL:98,FL:00}}). One idea could be to consider the distance
with respect to the operator norm of the difference of the orthogonal
projections $P_i$ onto the ranges of the frame transforms of two given frames
$\{ x_i \}_i$ and $\{ y_i \}_i$. Unfortunately, two frames of a certain
Hilbert space $H$ are similar if and only if these projections coincide
(Th.~7.2\cite{FL:98}). So we would only characterize classes of similar frames.
The better idea is to consider the operator norm distance of the respective
frame transforms $T_x$, $T_y$ or of the respective frame operators $S_x$,
$S_y$. The latter act as positive diagonalizable operators on the given
Hilbert C*-module, whereas the former map it to the standard countably
generated Hilbert C*-module $l_2(A)$.

For any tight frame $\{ y_i \}_i$ of a Hilbert space $H$ the corresponding
frame operator $S_y$ equals to the identity operator times the frame bound
value. For a given frame $\{ x_i \}_i$ of $H$ with frame bounds $C,D$ the
closest positive multiple of the identity operator to the frame operator
$S_x$ is $(1/C+1/D)/2 \cdot {\rm id}$. Unfortunately,
every tight frame $\{ y_i \}_i$ of $H$ with frame bound $(1/C+1/D)/2$ fulfils
this condition, and the relative position of tight frames in $H$ is of
greater importance than one can express that way. However, we got a hint for
the kind of factor to be used.
Now, let us consider the norm of differences of frame transforms.

\begin{proposition}
   Let $\{ x_i \}_i$ be a frame of a certain Hilbert space $H$ with frame bounds
   $C,D$, and let $S$ be its frame operator.
   The frame transform $T$ of the tight frame $\{ (\sqrt{C}+\sqrt{D})/2
   \cdot S^{1/2}(x_i) \}$ is the closest one w.r.t.~the operator norm to the
   frame transform $T_x$ of the given frame $\{ x_i \}_i$ among all the frame
   transforms of positive multiples of the normalized tight frame
   $\{ S^{1/2}(x_i) \}_i$.
\end{proposition}

\begin{proof}
Denote the frame transform of the tight frame $\{ \lambda \cdot S^{1/2}(x_i)
\}_i$ by $T_\lambda$ for $\lambda > 0$. Because the frame transforms
$T_\lambda$ and $T$ all possess the same coisometries in their respective
polar decompositions we obtain the equality $\| T_\lambda - T \| = \| \lambda
\cdot {\rm id} - S^{-1/2} \|$. Since the inequality $\sqrt{C} \leq S^{-1/2}
\leq \sqrt{D}$ is valid and both $\lambda \cdot {\rm id}$ and $S^{-1/2}$ are
diagonalizable with a common set of eigenvectors that forms a basis of $H$ we
obtain
  \[
   \| \lambda \cdot {\rm id} - S^{-1/2} \| = {\rm max} \, \{ |\lambda - \mu_j
   |\, : \, \mu_j \,\,\, {\rm any} \,\, {\rm eigenvalue} \,\, {\rm of} \,\,\,
   S^{-1/2} \} \, .
  \]
The right expression is minimal if and only if $\lambda$ is the arithmetic
mean of the lower and the upper spectral bound $\sqrt{C}$ and $\sqrt{D}$
of the positive invertible operator $S^{-1/2}$.
\end{proof}

Suppose the frame $\{ x_i \}_i$ of a certain Hilbert space $H$ is fixed and
possesses the frame bounds $C$ and $D$. Let us consider tight frames $\{ y_i
\}_i$ of $H$ with frame bound $(\sqrt{C}+\sqrt{D})/2$ for which the norm of
the difference of their frame transform $T_y$ and of the frame transform $T$
of the tight frame $\{(\sqrt{C}+\sqrt{D})/2 \cdot S_x^{1/2}(x_i) \}_i$ is
small. The next example shows that there are usually a lot of quite different
tight frames of $H$ with frame bound $(\sqrt{C}+\sqrt{D})/2$ the frame transforms
of which realize the same norm distance to the frame transform $T_x$ of the
initial frame as the distinguished tight frame $\{(\sqrt{C}+\sqrt{D})/2 \cdot
S_x^{1/2}(x_i) \}_i$ of $H$. We want to point out that the example works
already for finite-dimensional Hilbert spaces and for frames with finitely
many elements.

\begin{example} {\rm
  Let $H=l_2$ be a separable Hilbert space and fix its standard orthonormal
  basis $\{ e_i \}_i$. Set $x_1=e_1$, $x_2=3 \cdot e_2$ and $x_i = 2 \cdot e_i$
  for $i \geq 3$. The resulting set $\{ x_i \}_i$ is a Riesz basis of $H$ with
  frame bounds $C=1$ and $D=9$. Since $(\sqrt{C}+\sqrt{D})/2 = 2$ we can
  consider other tight frames $\{ y_i \}_i$ of $H$ with the same frame
  bound $4$ that are defined as $y_i = 2 \cdot e_i$ for $i \not= 3$ and
  $y_3 = 2{\rm e}^{{\rm i} \phi} \cdot e_3$ for some $\phi \in (-2 \cdot
  \arcsin (1/4), 2 \cdot \arcsin (1/4))$.
  Obviously, $\|T_x -T_y\| = \| S^{-1/2} - U \|$ for $T_x = VS_x^{-1/2}$ and
  $U=V^*T_y$, where $U$ maps $y_i$ to $\langle y_i,y_i \rangle e_i$, $i \in
  \mathbb N$ and $V$ is the identity map. Since both $S^{-1/2}$ and $U$ are
  normal, diagonalizable and commuting, with a basis consisting of common
  eigenvectors $\{ e_i \}_i$, we can estimate 
  \[
      \| T_x -T_y \| = \max \{ | \lambda_j - \mu_j | : j \in \mathbb N \, ,
      \,\, S^{-1/2}(e_j) = \lambda_j e_j \, , \,\, U(e_j) = \mu_j e_j \}
  \]
  see E.~A.~Azoff and C.~Davis\cite{AD}, or K.~R.~Davidson\cite{D1,D2}.
  The eigenvalues can be counted, they are $\lambda_i = \langle x_i,x_i
  \rangle^{1/2}$ and $\mu_i = \langle y_i,y_i \rangle^{1/2}$ for $i \not= 3$
  and $\mu_3= {\rm e}^{-{\rm i} \phi} \langle y_3,y_3 \rangle^{1/2}$.
  Taking the concrete values from the definitions of both the
  frames we obtain that the maximum of the difference of the corresponding
  eigenvalues is determined by the first two terms as long as $|\phi| < 2 \cdot
  \arcsin (1/4)$. So all these tight frames $\{ y_i \}_i$ parametrized by
  $\phi \in (-2 \cdot \arcsin (1/4), 2 \cdot \arcsin (1/4) )$ realize the same
  norm $\| T_x -T_y \| = 1$ for the difference of the respective frame
  transforms.}
\end{example}

Summarizing, the measure of nearness of a frame to some tight frame derived
from the norm of the difference of their frame transforms in general gives
an entire manifold of tight frames that are `closest' to a given frame with
respect to this measure. Moreover, if we allow the tight frame to span a
probably smaller Hilbert space than the original frame $\{ x_i \}_i$, then the
condition $D < 9/4 \cdot C$ to the frame bounds $C,D$ of $\{ x_i \}_i$
turns out to be occasionally essential to guarantee that the closest
tight frame spans exactly the same Hilbert space than the initial frame. In
other words, the distance between the square root of the lower frame bound $C$
and the arithmetic mean of the square roots of the lower and the upper frame
bounds $C$ and $D$, respectively, has to be smaller than the distance of the
square root of $C$ to zero.

The problem stated in the beginning of the present section remains unsolved
despite of the encouraging partial results indicated above,
even for the approximation of frames of Hilbert spaces by (normalized) tight
ones. We will continue our work to find a solution for it.

%%%%%%%%%%%%%%%%%%%%%%%%%%%%%%%%%%%%%%%%%%%%%%%%%%%%%%%%%%%%%%%%%%%%%%%%%%%%%
\acknowledgements

We are grateful to Radu Balan, Peter G.~Casazza, Mark Lammers, Christian
le Merdy and Vern I.~Paulsen for fruitful discussions and useful remarks
on the subjects explained in the present survey, especially on the use of
the modular point of view in classical frame theory.

%%%%%%%%%%%%%%%%%%%%%%%%%%%%%%%%%%%%%%%%%%%%%%%%%%%%%%%%%%%%%%%%%%%%%%%%%%%%%


\begin{thebibliography}{99}
\bibitem{FL:98} {M.~Frank, D.~R.~Larson}, ``Frames in Hilbert C*-modules and
  C*-algebras'', preprint, University of Houston, Houston, and Texas A{\&}M
  University, College Station, Texas, U.S.A., 1998.
\bibitem{FL:00} {M.~Frank, D.~R.~Larson}, ``A module frame concept for
  Hilbert C*-modules'', in: {\em Functional and Harmonic Analysis of Wavelets
  (San Antonio, TX, Jan.~1999)}, D.~R.~Larson, L.~W.~Baggett, eds.,
  AMS, Providence, R.I., {\em Contemp.~Math.} {\bf 247}, 207-233, 2000.
\bibitem{Wata} {Y.~Watatani}, ``Index for C*-subalgebras'', {\it Memoirs 
   Amer.~Math.~Soc.} {\bf 424}, 1990, 83 pp.
\bibitem{PP} {M.~Pimsner, S.~Popa}, ``Entropy and index for subfactors'',
  {\it Ann.~Scient.~Ec. Norm.~Sup.} {\bf 19}, 57-106, 1986.
\bibitem{BDH} {M.~Baillet, Y.~Denizeau, J.-F.~Havet}, ``Indice d'une
   esperance conditionelle'', {\it Compos.~Math.} {\bf 66}, 199-236, 1988.
\bibitem{FK} {M.~Frank, E.~Kirchberg}, ``On conditional expectations of
   finite index'', {\it J.~Oper. Theory} {\bf 40}, 87-111, 1998.
\bibitem{DPZ} {S.~Doplicher, C.~Pinzari, R.~Zuccante}, ``The C*-algebra of a
  Hilbert bimodule'', {\it Bollettino Unione Mat.~Ital., Sez.~B Artic.~Ric.~Mat.
  (8)} {\bf 1}, no.~2, 263-282, 1998.
\bibitem{KPW} {T.~Kajiwara, C.~Pinzari, Y.~Watatani}, ``Ideal structure and
  simplicity of the C*-algebras generated by Hilbert bimodules'', {\it
  J.~Funct.~Anal.} {\bf 159}, 295-322, 1998.
\bibitem{Badea} {C.~Badea}, ``The stable rank of topological algebras and a
  problem of R.~G.~Swan'', {\it J.~Funct.~Anal.} {\bf 160}, 42-78, 1998. 
\bibitem{Vill} {J.~Villadsen}, ``On the stable rank of simple C*-algebras'',
  {\it J.~Amer.~Math.~Soc.} {\bf 12}, 1091-1102, 1999.
\bibitem{CM} {A.~L.~Carey, V.~Mathai}, ``$L^2$-torsion invariants'', {\it
  J.~Funct.~Anal.} {\bf 110}, 377-409, 1992.
\bibitem{Lue} {W.~L\"uck}, ``Hilbert modules and modules over finite
  von Neumann algebras, and applications to $L^2$-invariants'',
  {\it Math.~Ann.} {\bf 309}, 247-285, 1997.
\bibitem{Fr} {M.~Frank}, ``Hilbertian versus Hilbert W*-modules, $L^2$- and
  other invariants'', preprint no.~9/2000, ZHS-NTZ, Universit\"at Leipzig /
  math.OA/0003185 (in xxx.lanl.gov), 2000.
\bibitem{Lance} {E.~C.~Lance}, {\it Hilbert C*-modules -- a Toolkit for
  Operator Algebraists}, {London Mathematical Society Lecture Note Series}
  {v.~210}, Cambridge University Press, Cambridge, England, 1995.
\bibitem{NEWO} {N.~E.~Wegge-Olsen}, {\it K-theory and C*-algebras
   -- a friendly approach}, Oxford University Press, Oxford, England,
   1993.
\bibitem{RW} {I.~Raeburn, D.~P.~Williams}, {\it Morita Equivalence and Continuous
  Trace C*-Algebras}, Mathematical Surveys and Monographs, v.~60,
  Amer.~Math.~Soc., Providence, R.I., 1998.
\bibitem{Fr92} {M.~Frank}, ``Direct integrals and Hilbert W*-modules
  (russ./engl.)'', in: Problems in Algebra, Geometry and Discrete Mathematics,
  eds.: O. B. Lupanov, A. I. Kostrikin, Moscow State University,
  Dept.~Mech.~Math., Moscow, Russia, 1992, 162-177 / e-print funct-an/9312002
  at xxx.lanl.gov. 
\bibitem{Frank:99} {M.~Frank}, ``Geometrical aspects of Hilbert
  C*-modules'', {\it Positivity} {\bf 3}, 215-243, 1999.
\bibitem{CHL} {P.~G.~Casazza, Deguang Han, D.~R.~Larson}, ``Frames for
  Banach spaces'', in: {\em Functional and Harmonic Analysis of Wavelets
  (San Antonio, TX, Jan.~1999)}, D.~R.~Larson, L.~W.~Baggett, eds.,
  AMS, Providence, R.I., {\em Contemp.~Math.} {\bf 247}, 149-182, 2000.
\bibitem{HL} {Deguang Han, D.~R.~Larson}, ``Frames, bases and
  group representations'', {\it Memoirs Amer. Math.~Soc.}, 1998, to appear.
\bibitem{Pa} {W.~L.~Paschke}, ``Inner product modules over B*-algebras'',
  {\em Trans.~Amer.~Math.~Soc.} {\bf 182}, pp.~443-468, 1973.
\bibitem{KaRi} {R.~V.~Kadison, J.~R.~Ringrose}, {\em Fundamentals of the Theory
  of Operator Algebras, I-IV}, Academic Press, New York, 1983.
\bibitem{Ped} {G.~K.~Pedersen}, ``Factorization in C*-algebras'',
  {\em Expos.~Math.} {\bf 16}, pp.~145-156, 1983.
\bibitem{Balan} {R.~Balan}, ``Equivalence relations and distances between
  Hilbert frames'', {\em Proc.~Amer.~Math.~Soc.} {\bf 127}, 2353-2366, 1999.
\bibitem{Tib} {T.~R.~Tiballi}, ``Symmetric orthogonalization of vectors in
  Hilbert spaces'', {\em Ph.D.~Thesis}, University of Houston, Houston,
  Texas, USA, 1991.
\bibitem{FPT} {M.~Frank, V.~I.~Paulsen, T.~R.~Tiballi}, ``Symmetric
  approximation of frames and bases in Hilbert spaces'',
  {\em Trans.~Amer.~Math.~Soc.}, to appear.
\bibitem{BF} {R.~Balan, M.~Frank}, ``An open problem'', Concentration Week,
  Texas A{\&}M University, College Station, TX, July 14-19, 1999, contribution.
\bibitem{Loew} {P.-O.~L\"owdin}, ``On the nonorthogonality problem'',
  {\em Adv.~Quantum Chem.} {\bf 5}, 185-199, 1970.
\bibitem{AD} {E.~A.~Azoff, C.~Davis}, ``On distances between unitary orbits
  of self-adjoint operators'', {\em Acta Sci.~Math.} {\bf 47}, 419-439,
  1984.
\bibitem{D1} {K.~R.~Davidson}, ``The distance between unitary orbits of
  normal operators'', {\it Acta Sci.~Math.} {\bf 50}, 213-223, 1986.
\bibitem{D2} {K.~R.~Davidson}, ``Estimating the distance between unitary
  orbits'', {\it J.~Oper.~Theory} {\bf 20}, 21-40, 1988.
\end{thebibliography}
\end{document}